# THE HOUGH TRANSFORM ESTIMATOR[1]

By Alexander Goldenshluger and Assaf Zeevi

*Haifa University and Columbia University*

This article pursues a statistical study of the Hough transform, the celebrated computer vision algorithm used to detect the presence of lines in a noisy image. We first study asymptotic properties of the Hough transform estimator, whose objective is to find the line that "best" fits a set of planar points. In particular, we establish strong consistency and rates of convergence, and characterize the limiting distribution of the Hough transform estimator. While the convergence rates are seen to be slower than those found in some standard regression methods, the Hough transform estimator is shown to be more robust as measured by its breakdown point. We next study the Hough transform in the context of the problem of detecting multiple lines. This is addressed via the framework of *excess mass functionals* and modality testing. Throughout, several numerical examples help illustrate various properties of the estimator. Relations between the Hough transform and more mainstream statistical paradigms and methods are discussed as well.

**1. Introduction.** The Hough transform (HT), due to Hough (1959), is one of the most frequently used algorithms in image analysis and computer vision [see, e.g., Ritter and Wilson (1996) and the survey articles by Leavers (1993) and Stewart (1999)]. The algorithm is most often used to detect and estimate parameters of multiple lines that are present in a noisy image (typically the image is first edge-detected and the resulting data serve as input to the algorithm).

In the particular case where only one line is present, the algorithm shares the same objective as simple linear regression, namely, estimating the slope and intercept of the line. While inference using regression methods is well understood, the statistical properties of the HT approach have not been

Received July 2002; revised October 2003.

[1]Supported in part by the German/Israeli G.I.F. Research Grant No. 2042-1126.4/2001.

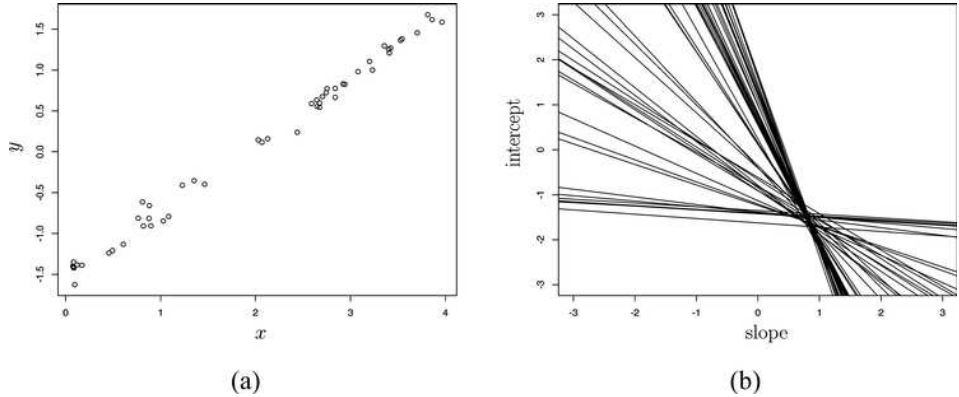

**(a)**                                              **(b)**

Fig. 1. *An illustration of the Hough transform:* (a) *the original scatterplot;* (b) *the Hough domain (dual plot).*

studied thoroughly. Most studies have focused almost exclusively on algorithmic and implementation aspects [for a comprehensive survey see, e.g., Leavers (1993)], while few articles pursue a statistical formulation [see, e.g., Kiryati and Bruckstein (1992) and Princen, Illingworth and Kittler (1994)].

The basic idea of the HT can be informally described as follows. Consider a set of planar points $\{(X_i, Y_i)\}_{i=1}^{n}$ depicted in Figure 1(a). The objective is to infer the parameters of the line that fits the data in the "best" manner. The key to the HT algorithm is to view each point as generating a line which is comprised of all pairs (slope, intercept) that are consistent with this point. Specifically, for the $i$th point this line is given by $L_i = \{(a, b) : Y_i = aX_i + b\}$. The set of random lines $\{L_i\}_{i=1}^{n}$ is plotted in the Hough domain, depicted in Figure 1(b). In the statistical literature this domain is referred to as the *dual plot*. Thus, co-linearity in the original set of points will manifest itself in a common intersection of lines in the dual plot.

In practice, the HT algorithm is implemented as follows. The Hough domain is first quantized into cells, and each such cell maintains a count of the number of lines that intersect it. The cell with the largest number of counts is the obvious estimator of the parameters of the original line. If one is focusing on detecting multiple lines, a threshold is specified and those cells with counts exceeding the threshold indicate the presence (and parametrization) of lines in the original image. A polar parametrization of the lines is also used in practical implementations, resulting in sinusoidal curves in the Hough domain [see, e.g., Ritter and Wilson (1996)].

The goal of this article is to provide analysis that formalizes and elucidates statistical properties of the HT methodology. The main contributions of this article are the following:

(i) We establish almost sure consistency of the HT estimator (Theorem 1), determine the rate of convergence and characterize the limiting

distribution (Theorem 2). The estimator is shown to have cube-root asymptotics [see, e.g., Kim and Pollard (1990)].

(ii) Robust properties of the HT estimator are derived. In particular, the breakdown point is determined (Theorem 3) and it is shown that this point can be made to be arbitrarily close to 50%. The theory is illustrated via a standard example.

(iii) We illustrate the effects of design parameters of the HT estimator on its performance via a simulation study.

(iv) We relate the multiple line detection problem to multi-modality testing in the Hough domain. In particular, asymptotic behavior of empirical excess mass functionals (Theorem 4) provides the building block by which one can pursue a test for the presence of multiple lines.

While a study focusing on the statistical properties of the HT is lacking in the literature, several strands of statistics-related research are akin to the HT approach. The concept of the dual plot has appeared already in early work of Daniels (1954), and in more recent work of Johnstone and Velleman (1985) and Rousseeuw and Hubert (1999). As we shall see in what follows, the HT estimator is closely related to regression methods such as least median of squares of Rousseeuw (1984), and $S$-estimators studied in Rousseeuw and Yohai (1984) and Davies (1990). Finally, the multiple line detection problem is intimately related to multi-modality testing using excess mass [see, e.g., Hartigan (1987), Müller and Sawitzki (1991) and Polonik (1995)]. The basic problem of estimating the location of a single mode studied by Chernoff (1964) can also be viewed as a one-dimensional application of the HT algorithm. Further details concerning some of these relations are given in the sequel.

The article has two main focal points: the first three sections, namely, Sections 2–4, focus on the HT estimator, while the subsequent Section 5 discusses testing of multiple lines. Section 2 describes the precise formulation of the HT estimator, while Section 3 studies large sample properties of the HT estimator (Section 3.1) and robustness (Section 3.2). Section 4 then focuses on some issues concerned with the design of the estimator, effects of the variates and relation of the method to other statistical approaches. The problem of testing for multiple lines is the subject of Section 5. Finally, Section 6 contains several concluding remarks. Proofs are collected in two appendices: Appendix A gives the proofs related to the properties of the HT estimator, while Appendix B contains the proofs related to the multiple line testing problem.

**2. Definition of the HT estimator.** Let data points $(X_1, Y_1), \ldots, (X_n, Y_n)$ be given on the plane. Each observation pair $(X_i, Y_i)$ defines a straight line in the Hough domain:

$$L_i : b = -X_i a + Y_i, \qquad i = 1, \ldots, n.$$

For a positive number $r$, let $B_r(\theta)$ denote the disc of radius $r$ centered at $\theta = (a, b)$. We are looking for a point $\hat{\theta} = (\hat{a}, \hat{b})$ in the Hough domain such that the maximal number of lines $L_i$ cross over the disc $B_r(\hat{\theta})$. More formally, the HT estimator $\hat{\theta}_{r,n}$ maximizes the objective function

$$M_{r,n}(\theta) := \frac{1}{n} \sum_{i=1}^{n} \mathbf{1}\{B_r(\theta) \cap L_i \neq \varnothing\}$$

with respect to $\theta = (a, b)$. Note that $L_i \cap B_r(\theta) \neq \varnothing$ if and only if the distance between the line $L_i$ and the disc center $\theta = (a, b)$ is less than or equal to $r$. Thus, $M_{r,n}(\theta)$ takes the following form:

$$(1) \qquad M_{r,n}(\theta) = \frac{1}{n} \sum_{i=1}^{n} \mathbf{1}\{|X_i a + b - Y_i|^2 \leq r^2(X_i^2 + 1)\},$$

and the HT estimator is defined by

$$(2) \qquad \hat{\theta}_{r,n} = \arg\max_{\theta \in \mathbb{R}^2} \frac{1}{n} \sum_{i=1}^{n} \mathbf{1}\{|X_i a + b - Y_i|^2 \leq r^2(X_i^2 + 1)\}.$$

Hence, $\hat{\theta}_{r,n}$ can be regarded as an $M$-estimator associated with the objective function $M_{r,n}(\cdot)$. Note that usually the above maximum is not unique; any point of the solution set may be chosen as $\hat{\theta}_{r,n}$. Note also that the above

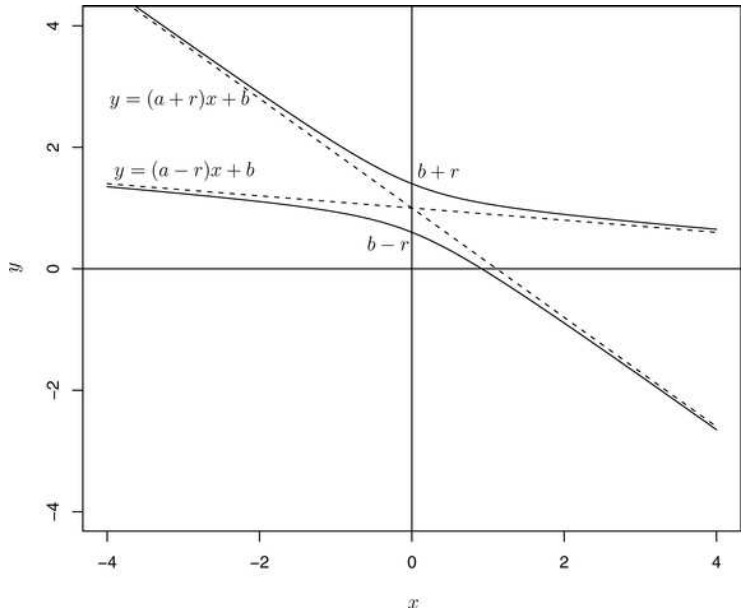

FIG. 2.   *Template of the HT estimator.*

definition of the HT estimator depends on the design parameter $r$. Denote by

$$(3) \qquad M_r(\theta) := \mathbb{E}M_{r,n}(\theta) = \mathbb{P}\{|Xa + b - Y|^2 \leq r^2(X^2 + 1)\}$$

the deterministic counterpart of $M_{r,n}(\theta)$.

The HT estimator admits the following geometrical interpretation. Let

$$(4) \qquad D_\theta = \{(x, y) : |xa + b - y|^2 \leq r^2(x^2 + 1)\}, \qquad \theta = (a, b) \in \mathbb{R}^2.$$

For given $\theta$, $D_\theta$ is the set of all points of the plane lying between two branches of a hyperbola that has straight lines $y = (a - r)x + b$ and $y = (a + r)x + b$ as its asymptotes; see Figure 2. Hence, the HT estimator given by (2) seeks the value $\theta$ such that the corresponding set $D_\theta$ covers the maximal number of data points. The set $D_\theta$ defines the so-called *template* of the HT in the observation space [e.g., Princen, Illingworth and Kittler (1992)]. We note that the template shape is determined by the choice of the *cell shape*, which is a disc of radius $r$ in our case. Various estimators may be defined using other cell shapes; the rectangular cell is most natural. However, the difference in properties of these estimators is marginal.

**3. Properties of the HT estimator.** Asymptotic properties of the HT estimator are studied under the following assumptions. Suppose that $(X_1, Y_1), \ldots, (X_n, Y_n)$ are independent identically distributed random observations drawn from the model

$$(5) \qquad\qquad Y = a_0 X + b_0 + \epsilon,$$

where:

(a) $X$ is independent of $\epsilon$, and

(b) $\epsilon$ is a random variable with bounded, symmetric and strictly unimodal density, $f(x) = f(-x) \ \forall x$.

By strict unimodality we mean that density $f$ has a maximum at a unique point, $x = 0$, and decreases in either direction as $x$ decreases or increases away from zero.

Let $\mathbb{P}_n$ denote the empirical measure of a sample of the pairs $(X_i, Y_i)$, $i = 1, \ldots, n$, and $\mathbb{P}$ be the common distribution of $(X_i, Y_i)$. Then the objective function $M_{r,n}(\theta)$ in (1) and its deterministic counterpart, $M_r(\theta)$, can be written as $M_{r,n}(\theta) = \mathbb{P}_n(D_\theta)$ and $M_r(\theta) = \mathbb{P}(D_\theta)$, where $D_\theta$ is defined by (4).

3.1. *Asymptotics.* We are interested in the asymptotic behavior of $\hat{\theta}_{r,n}$ as $n \to \infty$. The first theorem establishes consistency.

THEOREM 1. *Under assumptions* (a) *and* (b), *for any fixed* $r > 0$ *the estimator* $\hat{\theta}_{r,n}$ *is strongly consistent:*

$$\hat{\theta}_{r,n} \overset{a.s.}{\to} \theta_0 \qquad as \ n \to \infty, \ where \ \theta_0 = (a_0, b_0).$$

It is interesting to note that the consistency proof does not require existence of the expectation of the noise $\epsilon$. For example, the noise may be a sequence of i.i.d. Cauchy random variables. The next theorem establishes the asymptotic distribution of the centered and scaled estimator.

THEOREM 2. *Let $f$ be continuously differentiable with bounded first derivative, and let assumptions* (a) *and* (b) *hold. Assume that $X$ is a nondegenerate random variable with finite second moment. Then for every fixed $r > 0$,* $n^{1/3}(\hat{\theta}_{r,n} - \theta_0) \Rightarrow W$, *where $W$ has the distribution of the (almost surely unique) maximizer of the process* $\theta \mapsto \frac{1}{2}\theta^T V_0 \theta + G(\theta)$,

$$(6) \qquad V_0 = \mathbb{E}\{[f'(r\|Z\|) - f'(-r\|Z\|)]ZZ^T\},$$

$Z = (X, 1)^T$, *and $G$ is a zero-mean Gaussian process with continuous sample paths and stationary increments such that for any $\xi, \eta \in \mathbb{R}^2$,*

$$(7) \qquad \mathbb{E}[G(\xi) - G(\eta)]^2 = 2\mathbb{E}\{f(r\|Z\|)|Z^T(\xi - \eta)|\}.$$

The cube-root rates of convergence are due to the discontinuous nature of the objective function $M_{r,n}(\cdot)$. The most general results dealing with this type of asymptotics are given in Kim and Pollard (1990); see also van der Vaart and Wellner [(1996), Chapter 3]. Clearly the asymptotic distribution above is quite complicated. The one-dimensional instance, where $G(\cdot)$ is a Brownian motion, was first studied in Chernoff (1964) [see also, Groeneboom (1989) and Groenenboom and Wellner (2001)].

3.2. *Robustness.* One way to characterize the robustness of an estimator is through its breakdown properties. Intuitively, the *breakdown point* is the smallest amount of "contamination" necessary to "upset" an estimator entirely. We use the formal definition of the finite-sample breakdown point given by Donoho and Huber (1983). Let $\mathcal{Y}_n = \{(X_1, Y_1), \ldots, (X_n, Y_n)\}$ and $\hat{\theta} = \hat{\theta}(\mathcal{Y}_n)$ be an estimator based on $\mathcal{Y}_n$. Consider an additional data set $\mathcal{Y}'_k$ of size $k$. If by choice of $\mathcal{Y}'_k$ one can make $\hat{\theta}(\mathcal{Y}_n \cup \mathcal{Y}'_k) - \hat{\theta}(\mathcal{Y}_n)$ arbitrarily large, we say that $\hat{\theta}$ *breaks down* under contamination fraction $k/(n+k)$. The finite-sample *addition breakdown point* $\varepsilon_{\text{add}}(\hat{\theta}; \mathcal{Y}_n)$ is the minimal contamination fraction under which $\hat{\theta}$ breaks down:

$$\varepsilon_{\text{add}}(\hat{\theta}; \mathcal{Y}_n) = \min\left\{\frac{k}{n+k} : \sup_{\mathcal{Y}'_k} \|\hat{\theta}(\mathcal{Y}_n \cup \mathcal{Y}'_k) - \hat{\theta}(\mathcal{Y}_n)\| = \infty\right\}.$$

Similarly, the finite-sample *replacement breakdown point* of $\hat{\theta}$ is defined by

$$\varepsilon_{\mathrm{rep}}(\hat{\theta}; \mathcal{Y}_n) = \min\left\{\frac{k}{n} : \sup_{\mathcal{Y}_n^k} \|\hat{\theta}(\mathcal{Y}_n^k) - \hat{\theta}(\mathcal{Y}_n)\| = \infty\right\},$$

where $\mathcal{Y}_n^k$ denotes the corrupted sample obtained from $\mathcal{Y}_n$ by replacing $k$ data points of $\mathcal{Y}_n$ with arbitrary values. The following theorem gives the breakdown properties of the HT estimator $\hat{\theta}_{r,n}$.

THEOREM 3. *Let $\mathcal{Y}_n = \{(X_1, Y_1), \ldots, (X_n, Y_n)\}$ be a sample with no repeated values of $X$. Then*

$$\varepsilon_{\mathrm{add}}(\hat{\theta}; \mathcal{Y}_n) = \frac{\lfloor n M_{r,n}(\hat{\theta}_{r,n})\rfloor - 1}{n + \lfloor n M_{r,n}(\hat{\theta}_{r,n})\rfloor - 1},$$

$$\varepsilon_{\mathrm{rep}}(\hat{\theta}; \mathcal{Y}_n) = \frac{1}{n}\left\lfloor \frac{n M_{r,n}(\hat{\theta}_{r,n})}{2}\right\rfloor.$$

*Moreover, if the conditions of Theorem 1 hold, and the distribution of $X$ is continuous, then, as $n \to \infty$,*

$$\varepsilon_{\mathrm{add}}(\hat{\theta}_{r,n}; \mathcal{Y}_n) \overset{a.s.}{\to} p(1+p)^{-1}, \qquad \varepsilon_{\mathrm{rep}}(\hat{\theta}_{r,n}; \mathcal{Y}_n) \overset{a.s.}{\to} p/2,$$

*where $p = \mathbb{P}\{\epsilon^2 \leq r^2 \|Z\|^2\}$.*

We now turn to several remarks concerning the theorem. First, the assumption that the sample $\mathcal{Y}_n$ does not contain repeated observations of $X$ rules out parallel lines in the Hough domain. This assumption is quite typical in the context of the regression methods utilizing the dual plot approach [see, e.g., Daniels (1954)]. Second, the value of $r$ controls breakdown properties of the HT estimator: the larger $r$, the closer the breakdown point is to $1/2$. For example, if $r$ is chosen to be the $(1 - \alpha)$-quantile of the distribution of $\epsilon^2 \|Z\|^{-2}$, the addition breakdown point of the corresponding estimate is $(1 - \alpha)/(2 - \alpha)$ and the replacement breakdown point is $(1 - \alpha)/2$.

To illustrate the breakdown properties of the HT estimator, we consider a numerical example given in Rousseeuw (1984). The sample containing 30 "good" observations is generated from the model $Y_i = X_i + 2 + \epsilon_i$, where $\epsilon_i$ are Gaussian random variables with zero mean and standard deviation 0.2, and $X_i$ are uniformly distributed on $[1, 4]$. Then a cluster of 20 "bad" observations is added. These observations follow a bivariate Gaussian distribution with expectation $(7, 2)$ and covariance matrix $0.25I$. Figure 3 displays the data along with the least squares (LS), least median of squares (LMS) and the HT estimates. The LMS estimator is defined as the value of the parameter $\theta = (a, b)$ that minimizes the median$_{1 \leq i \leq n} |Y_i - aX_i - b|^2$ [see Rousseeuw (1984)]. The parameter $r$ of the HT estimator is set to 0.15. Under conditions of the experiment $\mathbb{P}\{\epsilon^2(X^2 + 1)^{-1} \leq 0.15^2\} \approx 0.923$, which approximately corresponds to a 46% replacement breakdown point. The HT estimator is calculated by direct maximization of (2) on the square $[-3, 3] \times [-3, 3]$ using a uniform rectangular grid comprised of 250,000 points. Because the solution is not unique, the average of the grid points where the maximum is achieved is taken as the estimate. Thus, the HT estimate yields $\hat{a} = 0.917$ and $\hat{b} = 2.173$, which is quite close to the original values $a_0 = 1$ and $b_0 = 2$. In general, behavior of the HT estimate in this example is very similar to that of the LMS.

## 4. Discussion.

4.1. *Choice of the radius* $r$. The properties of the HT estimator depend on the choice of a parameter $r$. The results of Section 3 assert that the maximum HT estimator is consistent for any choice of $r$, and the asymptotic distribution is given in Theorem 2. Thus, a reasonable choice of $r$ would be the value minimizing the variance of the limiting random variable in Theorem 2. Unfortunately, the asymptotic distribution is not tractable, and we cannot use it as a basis to make a choice of $r$. Clearly, large values of $r$ lead to a large connected solution set, and in this case the estimation accuracy depends crucially on the way the estimator is chosen from the solution set. On the other hand, small values of $r$ lead to an "under-smoothed" dual plot, and the solution set is a union of many disconnected sets. In this case estimation accuracy of the average estimator may be very poor.

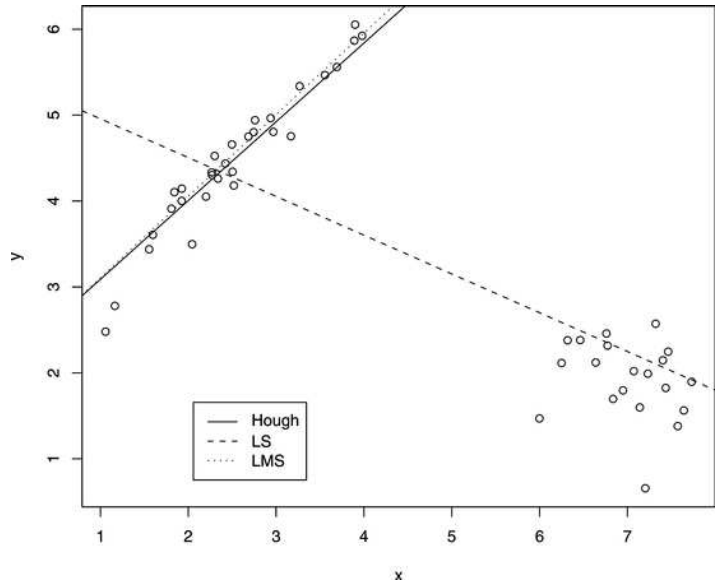

FIG. 3.  *An illustration of the breakdown properties of the HT estimator. The data set consists of* 30 *observations from the underlying linear regression model and* 20 *"bad" data points (the cluster on the right).*

To study how estimation accuracy depends on $r$, we conducted the following simulation experiment. For sample sizes $n = 25, 50, 100$ we generate data sets from the model $Y_i = X_i + 2 + \epsilon_i$, where $\epsilon_i$ are Gaussian random variables with zero mean and standard deviation 0.5, and $X_i$ are uniformly distributed on $[-2, 2]$. The HT estimator is computed for different values of $r$. In our implementation we used the square $[-3, 3] \times [-3, 3]$ as the search region. The value of the objective function is computed at nodes of the regular grid comprised of 360,000 points. The resulting HT estimator is set to be the average of the grid nodes where the maximum of the objective function is achieved. Simulation results are given in Table 1. The table presents the values of the HT estimates of the parameters $(a_0, b_0) = (1, 2)$ averaged over 1,000 replications, along with the square root of the resulting mean squared error. Closer inspection of the results shows that the mean squared error first decreases as $r$ grows, but when $r$ becomes large, an increase in the mean squared error is observed. Overall, it seems that the estimation accuracy is relatively stable as $r$ varies over a wide range of values. This phenomenon has been consistently observed for various data sets generated from different models. (The results described in Table 1 are one such representative example.) Finally, we note that, in practice, it may be advantageous to take $r$ slowly tending to zero as $n \to \infty$. This might be particularly important in the problem of multiple line testing discussed in Section 5. However, analysis of theoretical properties of such an estimator is beyond the scope of this article.

4.2. *Equivariance properties and the effect of design variables.* We now briefly mention some equivariance properties of the HT estimator. In the context of regression estimators, different notions of equivariance are considered [see, e.g., Rousseeuw and Leroy (1987), page 116]. An estimator $\hat{\theta}$ is said to be *regression equivariant* if $\hat{\theta}(\{X_i, Y_i + cX_i\}_{i=1}^n) = \hat{\theta}(\{X_i, Y_i\}_{i=1}^n) + c$, where $c$ is an arbitrary constant. It is *scale equivariant* if $\hat{\theta}(\{X_i, cY_i\}_{i=1}^n) = c\hat{\theta}(\{X_i, Y_i\}_{i=1}^n)$ and *affine equivariant* if $\hat{\theta}(\{cX_i, Y_i\}_{i=1}^n) = c^{-1}\hat{\theta}(\{X_i, Y_i\}_{i=1}^n)$ for $c \neq 0$.

It is easily seen that the HT estimator $\hat{\theta}_{r,n}$ is *regression equivariant*, but not *scale* and *affine equivariant*. The equivariance properties of the HT estimator are clearly intimately related to the Hough template. In particular, the template displayed in Figure 2 implies that the estimate treats differently observations with small and large $X$-variate values. The straight lines

TABLE 1

*Estimation accuracy of the HT estimator. The numbers in parenthesis are the (slope, intercept) estimates, and the value below them is the associated root mean squared error. All values are obtained by averaging over 1000 replications*

| | Sample size | | |
|---|---|---|---|
| $r$ | $n = 25$ | $n = 50$ | $n = 100$ |
| 0.025 | (0.992, 1.981) | (0.995, 2.009) | (0.990, 2.009) |
| | 0.407 | 0.297 | 0.245 |
| 0.04 | (0.999, 1.989) | (0.997, 2.013) | (0.995, 2.013) |
| | 0.392 | 0.284 | 0.231 |
| 0.05 | (1.003, 2.001) | (0.995, 2.001) | (1.001, 2.009) |
| | 0.354 | 0.272 | 0.219 |
| 0.075 | (1.011, 2.007) | (0.992, 2.008) | (1.000, 2.009) |
| | 0.322 | 0.264 | 0.213 |
| 0.1 | (1.009, 2.008) | (0.996, 2.009) | (0.998, 2.015) |
| | 0.308 | 0.251 | 0.204 |
| 0.2 | (1.000, 2.010) | (0.997, 2.012) | (1.000, 2.004) |
| | 0.264 | 0.208 | 0.164 |
| 0.4 | (1.001, 2.010) | (0.999, 2.007) | (0.996, 2.003) |
| | 0.220 | 0.171 | 0.137 |
| 0.5 | (0.996, 2.008) | (0.995, 2.004) | (0.994, 2.001) |
| | 0.211 | 0.174 | 0.135 |
| 0.75 | (1.012, 1.999) | (1.002, 1.996) | (0.999, 2.003) |
| | 0.248 | 0.209 | 0.172 |
| 0.8 | (1.015, 1.997) | (1.002, 1.997) | (0.996, 2.002) |
| | 0.254 | 0.219 | 0.179 |

in the Hough domain corresponding to the observations with large $X_i$ values are very steep. If the majority of the observations have a large $X$-coordinate and the standard deviation of the noise is small, then the corresponding straight lines are nearly parallel. In this case behavior of the HT estimator may be quite poor.

To illustrate the effect of the design distribution, we generate 100 independent observations from the model $Y_i = X_i + 2 + \epsilon_i$, where $\epsilon_i$ are Gaussian random variables with zero mean and standard deviation 0.5. Figure 4 displays the perspective plots of the objective function $M_{0.3,n}(\theta)$, along with the corresponding dual plots for two different design distributions. Figure 4(a) and (b) corresponds to the explanatory variables $X_i$ uniformly distributed on $[-2, 2]$, while Figure 4(c) and (d) shows the case of $X_i$ uniformly distributed on $[20, 24]$. In the second case the objective function is very flat. This leads to a large solution set and high variability of the HT estimator. Theoretically, when $X_i$ are large, the matrix $V_0$ appearing in (6) is nearly singular because $f'(r\|Z\|) - f'(-r\|Z\|)$ is close to zero. Therefore, the asymptotic distribution of $\hat{\theta}_{r,n}$ is close to the distribution of the point of maximum of a zero mean Gaussian process given in (7). To recapitulate this point, the influence of the design distribution on estimation accuracy suggests that it would be reasonable, in practice, to center the explanatory variables before applying the HT estimator. We note that in computer vision applications this does not typically pose a problem as the measurement units used for the $X$-coordinate are image-independent.

4.3. *Related regression methods.* The HT estimator may be viewed as a counterpart to an $S$-estimator [cf. Rousseeuw and Yohai (1984) and Davies (1990)]. Indeed, fix $\delta \in (0,1)$ and consider the following optimization prob-

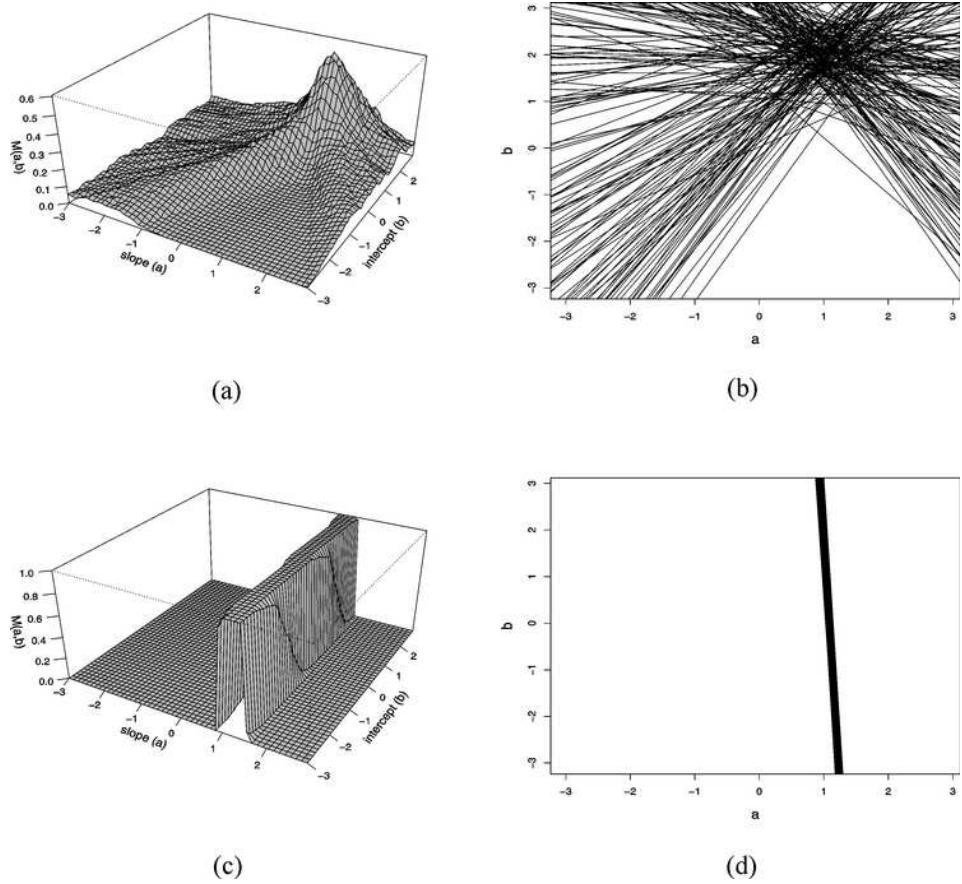

FIG. 4.  *Perspective plots of $M_{r,n}(\theta)$ along with the corresponding dual plots:* (a), (b) $X_i$ *are uniformly distributed on* $[-2, 2]$; (c), (d) $X_i$ *are uniformly distributed on* $[20, 24]$.

lem:

$$\mathcal{P}(\delta): \qquad \min_{\theta=(a,b)\in\mathbb{R}^2} r$$

(8)

$$\text{s.t.} \qquad M_{r,n}(\theta) = \frac{1}{n}\sum_{i=1}^{n}\mathbf{1}\{|Y_i - aX_i - b|^2 \le r^2(X_i^2+1)\} \ge 1-\delta.$$

Solution of (8) defines the $S$-estimator $\hat{\theta}_{\delta,n}$ whose replacement breakdown point equals $\varepsilon_{\mathrm{rep}}(\hat{\theta}_{\delta,n};\mathcal{Y}_n) = \min(\delta, 1-\delta)$ [cf. Davies (1990)]. The LMS estimator, see Rousseeuw (1984), can be written in a form similar to (8). In this specific case $\delta = n^{-1}(\lfloor n/2 \rfloor + 1)$ and $X_i^2 + 1$ on the right-hand side should be replaced by 1. Recall that, by definition, the HT estimator $\hat{\theta}_{r,n}$ solves the following optimization problem:

$$\mathcal{Q}(r): \max_{\theta=(a,b)\in\mathbb{R}^2} M_{r,n}(\theta) = \frac{1}{n}\sum_{i=1}^{n}\mathbf{1}\{|Y_i - aX_i - b|^2 \le r^2(X_i^2+1)\}.$$

Then the connection between the HT estimator and the $S$-estimator (8) is as follows. For a given $\delta > 0$, let $\hat{r} = \mathrm{val}(\mathcal{P}(\delta))$, where $\mathrm{val}(\cdot)$ is the value of the optimization problem, and let $\hat{\theta}_{\delta,n}$ be the solution to $\mathcal{P}(\delta)$. Then, clearly $\mathrm{val}(\mathcal{Q}(\hat{r})) \ge 1-\delta$, and $\hat{\theta}_{\hat{r},n}$ belongs to the solution set of $\mathcal{P}(\delta)$. Thus, with this particular choice of $r$, the HT estimator and the corresponding $S$-estimator are identical; in particular, $\varepsilon_{\mathrm{rep}}(\hat{\theta}_{\hat{r},n};\mathcal{Y}_n) = \min(\delta, 1-\delta)$.

**5. Multiple line detection.** In practice, the Hough domain is discretized into cells, and the number of lines crossing each cell is counted. Next, each of the cells is examined to search for "high counts." In particular, cells with counts exceeding some predetermined threshold correspond to "detected" lines in the original space. This procedure amounts to an exhaustive search for local maxima (threshold crossings) in the Hough domain. Thus, in contrast to other line fitting procedures, the HT is used to estimate several lines simultaneously. It should be noted, however, that points of local maxima do not necessarily correspond to actual line parameters. Consequently, in the case of multiple lines it is more accurate to view the HT as a tool for testing or *detecting* the presence of straight lines in images. This has also been recognized in the computer vision literature [cf. Princen, Illingworth and Kittler (1994)].

In view of the above, one can view the multiple line detection problem using the HT as testing for multi-modality in the Hough domain. Testing multi-modality is a subject of vast literature. This problem is characterized by the fact that only one-sided inference is possible [see, e.g., Donoho (1988)], that is, the only verifiable hypotheses are of the type "there are at least three lines in the image." The most appropriate approach for our purposes is based on the concept of *excess mass* [see Hartigan (1987), Müller and Sawitzki (1991) and Polonik (1995)], which is typically used in the "mode testing" problem. In the context of the HT, this excess mass corresponds to regions in the parameter space (Hough domain) where large counts are present.

5.1. *Excess mass functionals.* Let $(X_1, Y_1), \ldots, (X_n, Y_n)$ be a sample of i.i.d. random variables, and, for $r > 0$ and $\theta = (a, b) \in \mathbb{R}^2$, let $M_{r,n}(\theta)$ and $M_r(\theta)$ be given by (1) and (3), respectively. We stress that $(X_1, Y_1), \ldots, (X_n, Y_n)$ are not assumed to be drawn from the linear model (5). Throughout this section we suppose that parameter $\theta$ is confined to a compact set $\Theta_0 \subset \mathbb{R}^2$.

*The excess mass functional* is defined by

$$E(\lambda) := \int (M_r(\theta) - \lambda)^+ \, d\theta$$

$$= \int_{\Theta_\lambda} M_r(\theta) \, d\theta - \lambda \mathcal{L}\{\Theta_\lambda\},$$

where $(x)^+ := \max(0, x)$, $\Theta_\lambda := \{\theta \in \mathbb{R}^2 : M_r(\theta) \geq \lambda\}$, and $\mathcal{L}\{\cdot\}$ stands for Lebesgue measure in $\mathbb{R}^2$. We call $\Theta_\lambda$ the $\lambda$-*level set*; note that $\Theta_\lambda$ is closed and bounded because $M_r(\cdot)$ is continuous. For a compact set $\Theta \subset \mathbb{R}^2$ and $\lambda \in (0, 1)$, let us define

$$H_\lambda\{\Theta\} := \int_\Theta M_r(\theta) \, d\theta - \lambda \mathcal{L}\{\Theta\}.$$

Then $E(\lambda) = \sup\{H_\lambda\{\Theta\} : \Theta \subset \mathbb{R}^2 \text{ compact}\}$. The empirical version of the excess mass functional is obtained by substituting $M_{r,n}(\cdot)$ for $M_r(\cdot)$ in the definition, namely,

$$E_n(\lambda) := \int (M_{r,n}(\theta) - \lambda)^+ \, d\theta$$

$$= \int_{\Theta_{\lambda,n}} M_{r,n}(\theta) \, d\theta - \lambda \mathcal{L}\{\Theta_{\lambda,n}\},$$

where $\Theta_{\lambda,n} = \{\theta \in \mathbb{R}^2 : M_{r,n}(\theta) \geq \lambda\}$ is *the empirical $\lambda$-level set.* Using the notation

$$H_{\lambda,n}\{\Theta\} := \int_\Theta M_{r,n}(\theta) \, d\theta - \lambda \mathcal{L}\{\Theta\},$$

we have that $E_n(\lambda) = \sup\{H_{\lambda,n}\{\Theta\} : \Theta \subset \mathbb{R}^2, \text{ compact}\}$. Note that the empirical $\lambda$-level set $\Theta_{\lambda,n}$ is a closed subset of $\mathbb{R}^2$; this follows from the fact that $M_{r,n} : \mathbb{R}^2 \to [0, 1]$ is upper semi-continuous [see, e.g., Rudin (1987), pages 37 and 38]. Since the parameter $\theta$ is assumed to take values in the compact set $\Theta_0$, $\Theta_{\lambda,n}$ is also bounded.

Following Polonik ([1995](#)), we also consider the excess mass functional over some classes of subsets in $\mathbb{R}^2$. Let $\mathcal{T}$ be a class of compact subsets of $\mathbb{R}^2$. *The excess mass functional over $\mathcal{T}$ at level $\lambda \in (0,1)$ is given by*

$$E_\mathcal{T}(\lambda) := \sup\{H_\lambda\{\Theta\} : \Theta \in \mathcal{T}\} = \sup_{\Theta \in \mathcal{T}}\left[\int_\Theta M_r(\theta)\,d\theta - \lambda \mathcal{L}\{\Theta\}\right].$$

Every set $\Theta_\lambda(\mathcal{T}) \in \mathcal{T}$ satisfying $E_\mathcal{T}(\lambda) = H_\lambda\{\Theta_\lambda(\mathcal{T})\}$ is called *the $\lambda$-level set in $\mathcal{T}$*. Clearly, $E_\mathcal{T}(\lambda) \le E(\lambda)$ and $E_\mathcal{T}(\lambda) = E(\lambda)$ if $\Theta_\lambda \in \mathcal{T}$. The empirical version $E_{\mathcal{T},n}(\lambda)$ of $E_\mathcal{T}(\lambda)$ is defined by

$$E_{\mathcal{T},n}(\lambda) := \sup\{H_{\lambda,n}\{\Theta\} : \Theta \in \mathcal{T}\}$$

$$= \int_{\Theta_{\lambda,n}(\mathcal{T})} M_{r,n}(\theta)\,d\theta - \lambda \mathcal{L}\{\Theta_{\lambda,n}(\mathcal{T})\},$$

where $\Theta_{\lambda,n}(\mathcal{T})$ is the *empirical $\lambda$-level set in $\mathcal{T}$*.

We stress that the excess mass approach is very natural in the context of the HT. In particular, the value of $E_n(\lambda)$ conveniently quantifies the total sum of counts corresponding to cells with counts exceeding $\lambda$. Consequently, asymptotic behavior of the empirical excess mass functional is of interest.

5.2. *Asymptotics of the empirical excess mass functional.* The asymptotic behavior of the empirical excess mass functional is the key building block in a statistical procedure for detecting multiple lines; this is given in the next theorem. To that end, let us denote

$$\nu_n(\lambda) := \sqrt{n}\int_{\Theta_\lambda}[M_{r,n}(\theta) - M_r(\theta)]\,d\theta, \qquad \lambda \in \Lambda := [\underline{\lambda}, \overline{\lambda}] \subset (0,1),$$

and let $l^\infty(\Lambda)$ denote the space of all uniformly bounded real-valued functions over $\Lambda$.

THEOREM 4.    *Suppose that $M_r : \mathbb{R}^2 \mapsto [0,1]$ satisfies*

$$(9) \qquad \lim_{\delta \to 0}\sup_{\lambda \in \Lambda}\mathcal{L}\{\{\theta : |M_r(\theta) - \lambda| < \delta\}\} = 0.$$

*Then:*

(i) $\sup_{\lambda \in \Lambda}|\sqrt{n}[E_n(\lambda) - E(\lambda)] - \nu_n(\lambda)| = o_p(1)$ *as $n \to \infty$, and*

$$(10) \qquad \nu_n(\lambda) \Rightarrow \int_{\Theta_\lambda} G(\theta)\,d\theta \qquad in\ \ell^\infty(\Lambda), n \to \infty,$$

*where $G(\cdot)$ is a zero mean Gaussian random field with covariance kernel*

$$(11) \qquad \mathbb{E}[G(\xi)G(\eta)] = \mathbb{P}\{|Z^T\xi - Y| \le r\|Z\|, |Z^T\eta - Y| \le r\|Z\|\}$$

$$- \mathbb{P}\{|Z^T\xi - Y| \le r\|Z\|\}\mathbb{P}\{|Z^T\eta - Y| \le r\|Z\|\},$$

*where $Z = (X,1)^T$ and $\xi, \eta \in \mathbb{R}^2$.*

(ii) *Let $\mathcal{T}$ denote the class of compact subsets of $\mathbb{R}^2$ such that $\Theta_\lambda \in \mathcal{T}$ for every $\lambda \in \Lambda$. Then*

$$\sup_{\lambda \in \Lambda} |\sqrt{n}[E_{\mathcal{T},n}(\lambda) - E(\lambda)] - \nu_n(\lambda)| = o_p(1), \qquad n \to \infty,$$

*and* (10) *holds.*

The asymptotics of the empirical excess mass functional are determined by two factors: the asymptotic behavior of the random field $M_{r,n}(\theta)$ and the asymptotic behavior of the (random) level set $\Theta_{\lambda,n}$. There are essentially two main ideas that underlie the proof: (i) the class of sets generated by the Hough template, $\mathcal{D} = \{D_\theta : \theta \in \mathbb{R}^2\}$, is a separable VC class of sets, and, thus, a uniform central limit theorem holds for the random field $M_{r,n}(\cdot)$ [cf. Proposition 2]; (ii) under assumption (9), which essentially posits that the deterministic field $M_r(\cdot)$ does not have "flat parts," the convergence of the random field also implies convergence of the associated (random) level sets to their deterministic counterparts. In the absence of assumption (9), difficulties can easily arise in "mode testing" [see Müller and Sawitzki (1991) and Polonik (1995), where a similar condition is imposed in the context of excess mass testing for modes of a distribution].

5.3. *Testing for multiple lines.* We now sketch how Theorem 4 may be used for detecting multiple lines in some specific cases. To illustrate the ideas, consider the following hypothesis test:

(12)                $H_0$: one line    vs.    $H_1$: more than one line.

The rigorous interpretation of the above is that "under the null hypothesis," the data is generated by the model (5) with some unknown $\theta_0 = (a_0, b_0)$, and assumptions (a) and (b) of Section 3 hold. To characterize the behavior of excess mass functionals under the null hypothesis, we will need the next result which essentially states that under $H_0$ the $\lambda$-level set $\Theta_\lambda$ for $\lambda \in \Lambda$ is a convex set which is balanced around $\theta_0 = (a_0, b_0)$.

PROPOSITION 1. *Assume that the data are generated by the model* (5), *and assumptions* (a) *and* (b) *hold. Then $M_r(\theta) = \tilde{M}_r(\theta - \theta_0)$ for some function $\tilde{M}_r(\cdot)$ which is symmetric near zero with unique mode at $\theta = 0$. In addition, the set $\tilde{\Theta}_\lambda = \{\theta \in \mathbb{R}^2 : \tilde{M}_r(\theta) \geq \lambda\}$ is a closed convex and balanced set (i.e., if $\theta \in \tilde{\Theta}_\lambda$, then $-\theta \in \tilde{\Theta}_\lambda$).*

First consider the testing problem under the assumption that the distributions of $\epsilon$ and $X$ are known. Suppose that $M_r(\cdot)$ has no "flat parts," that is, (9) holds. By Proposition 1, under $H_0$ the excess mass functional $E(\lambda)$ is completely specified and given by

$$E_*(\lambda) = \int (M_r(\theta) - \lambda)^+ \, d\theta$$

$$= \int (\tilde{M}_r(\theta) - \lambda)^+ \, d\theta$$

$$= \int (\mathbb{P}\{|\epsilon + Z^T\theta| \le r\|Z\|\} - \lambda)^+ \, d\theta.$$

Thus, (12) reduces to testing

$$H_0' : E(\lambda) = E_*(\lambda) \qquad \forall \lambda \in \Lambda \quad \text{vs.} \quad H_1' : E(\lambda) \ne E_*(\lambda) \qquad \text{for some } \lambda \in \Lambda.$$

It follows from Theorem 4(i) that

$$T_n' := \sqrt{n} \sup_{\lambda \in \Lambda} |E_n(\lambda) - E_*(\lambda)| \Rightarrow \chi,$$

where $\chi := \sup_{\lambda \in \Lambda} |\int_{\Theta_\lambda} G(\theta) \, d\theta|$. Observe that $\Theta_\lambda = \tilde{\Theta}_\lambda + \theta_0$, hence,

$$\chi = \sup_{\lambda \in \Lambda} \left| \int_{\tilde{\Theta}_\lambda} G(\theta - \theta_0) \, d\theta \right| = \sup_{\lambda \in \Lambda} \left| \int_{\tilde{\Theta}_\lambda} \tilde{G}(\theta) \, d\theta \right|,$$

where $\tilde{G}(\cdot) = G(\cdot - \theta_0)$. We note that the covariance kernel of the zero mean Gaussian process $\tilde{G}(\cdot) := G(\cdot - \theta_0)$ does not depend on $\theta_0$ and is given by (11) with $Y$ replaced by $\epsilon$. Thus, the test can be based on the statistic $T_n'$ whose asymptotic distribution does not depend on unknown parameter $\theta_0$, and is completely specified under $H_0$, provided that the distributions of $\epsilon$ and $X$ are known. Such a test will be consistent against all alternatives of the type $|E(\lambda) - E_*(\lambda)| > 0$ for some $\lambda \in \Lambda$. We note that although the assumption that the distributions of $\epsilon$ and $X$ are known may seem to be restrictive, it is quite typical in many application settings [see, e.g., Princen, Illingworth and Kittler (1994)].

If the distributions of $X$ and $\epsilon$ are unknown, $T_n'$ cannot be computed and, therefore, testing the presence of one line against multiple lines is more complicated. In this setting one can pursue the multiple line testing problem by comparing restricted and unrestricted empirical excess mass functionals. Proposition 1 states that under the null hypothesis, the $\lambda$-level set $\Theta_\lambda$ is convex and balanced around $\theta_0$. Therefore, the test may be based on comparing $E_n(\lambda)$ with the empirical excess mass $E_{\mathcal{C},n}(\lambda)$ over the set $\mathcal{C}$ of all compact convex subsets of $\mathbb{R}^2$. Thus, we consider testing

$$\tilde{H}_0' : \Theta_\lambda \in \mathcal{C} \qquad \forall \lambda \in \Lambda \quad \text{vs.} \quad \tilde{H}_1' : \Theta_\lambda \notin \mathcal{C} \qquad \text{for some } \lambda \in \Lambda.$$

In view of Theorem 4, a natural test statistic is $\tilde{T}_n' := \sqrt{n} \sup_{\lambda \in \Lambda} |E_n(\lambda) - E_{\mathcal{C},n}(\lambda)|$, and $\tilde{H}_0'$ should be rejected for large values of $\tilde{T}_n'$. Under $\tilde{H}_0'$, $\tilde{T}_n' = O_p(1)$

as $n \to \infty$. On the other hand, if $E(\lambda) - E_{\mathcal{C}}(\lambda) > 0$ for some $\lambda \in \Lambda$, then by Theorem 4 the power of the test based on $\tilde{T}'_n$ converges to 1 as $n \to \infty$. Thus, the described test is consistent against all alternatives of the type $E(\lambda) - E_{\mathcal{C}}(\lambda) > 0$ for some $\lambda \in \Lambda$. Unfortunately, the limiting distribution of $\tilde{T}'_n$ is not available; in general, it depends on the rate at which $\sup_{\lambda \in \Lambda} \mathcal{L}\{\{\theta : |M_r(\theta) - \lambda| < \delta\}\}$ goes to zero as $\delta \to 0$ [cf. (9)]. We note that even though the condition $E(\lambda) - E_{\mathcal{C}}(\lambda) > 0$ does not imply that $\Theta_\lambda \neq \Theta_\lambda(\mathcal{C})$, in many situations this is the case.

## 6. Concluding remarks.

1. The HT estimator can be used in the multiple regression context. Assume the model

$$Y = \sum_{k=1}^{p} \beta_k X_k + \epsilon,$$

and denote $\theta = (\beta_1, \ldots, \beta_p)^T$ and $Z = (X_1, \ldots, X_p)^T$. Then the HT estimator is defined by

$$(13) \qquad \hat{\theta}_{r,n} = \arg\max_{\theta \in \mathbb{R}^p} \frac{1}{n} \sum_{i=1}^{n} \mathbf{1}\{|Y_i - \theta^T Z_i|^2 \leq r^2 \|Z_i\|^2\}.$$

It can be easily seen that Theorems 1–3 hold for the multiple regression setup with obvious modifications. In particular, the breakdown point given in Theorem 3 does not depend on the dimension. Unfortunately, the maximization problem in (13) is difficult and cannot be solved as easily as in the two-dimensional case.

2. The slow, cube root, convergence rate of the HT estimator is a consequence of the discontinuous objective function. Kim and Pollard (1990) study this phenomenon and survey various estimation settings in which cube root convergence rates govern the asymptotics. To this end, the original objective function might be approximated by a smooth function, and the resulting modified "smoothed" estimator would have standard $\sqrt{n}$ asymptotics and "good" breakdown properties. In this case maximization of the objective function can be pursued using a gradient-based search.

3. A variety of modified estimators may be obtained using different cell shapes in the Hough domain. For example, a vertical line segment of length $2r$ as a cell shape in the Hough domain corresponds to an estimator which maximizes

$$\frac{1}{n} \sum_{i=1}^{n} \mathbf{1}\{|Y_i - \theta^T Z_i|^2 \leq r^2\}$$

over $\theta \in \mathbb{R}^2$. The template of this estimator represents a strip of width $2r$ measured in the vertical direction. Such an estimator can be viewed as a counterpart to the LMS estimator. The properties of the estimator are quite similar to those of the HT estimator. In addition, such an estimator is scale and affine equivariant.

4. Fitting a straight line when both variables are subject to random errors can be treated using the described techniques. For example, it can be easily shown that the estimator based on the vertical line-segment cell is consistent, provided the errors have symmetric strongly unimodal densities.

## APPENDIX A: PROOFS FOR SECTION 3

PROOF OF THEOREM 1.  Conditioning on $X$, we have for $\theta \neq \theta_0$,

$$
\begin{aligned}
\mathbb{E}[M_{r,n}(\theta)|X] &= \mathbb{P}\{|Xa + b - Y|^2 \leq r^2(X^2 + 1)|X\} \\
&= \mathbb{P}\{-r\sqrt{X^2 + 1} - X(a - a_0) - (b - b_0) \\
&\qquad \leq -\epsilon \leq r\sqrt{X^2 + 1} - X(a - a_0) - (b - b_0)|X\} \\
&< \mathbb{P}\{-r\sqrt{X^2 + 1} \leq -\epsilon \leq r\sqrt{X^2 + 1}|X\}.
\end{aligned}
$$

The last inequality is a consequence of the Anderson lemma [Anderson (1955)] and the fact that $f$ is symmetric and strictly unimodal. Hence, $\theta_0$ is a unique point of maximum of function $M_r(\theta) := \mathbb{E}M_{r,n}(\theta)$ for any $r > 0$. In particular, denoting by $B_\varepsilon(\theta_0)$ the ball of radius $\varepsilon$ with center $\theta_0$, we have that for any $\varepsilon > 0$,

$$
\text{(14)} \qquad \max_{\theta \in B_\varepsilon^c(\theta_0)} M_r(\theta) < M_r(\theta_0).
$$

The point of maximum of $M_r(\cdot)$ is, thus, unique and well separated.

Consider the class of sets $\mathcal{D} = \{D_\theta, \theta \in \mathbb{R}^2\}$, where $D_\theta$ is defined in (4). This class has polynomial discrimination, that is, it is a Vapnik–Cervonenkis (VC) class of sets [see Pollard (1984), Definition II.13, or van der Vaart and Wellner (1996), page 85]. Indeed, as was mentioned before, $\mathcal{D}$ is a class of subsets of the plane generated by a linear space of quadratic forms. Hence, by Lemma II.18 in Pollard (1984), $\mathcal{D}$ has polynomial discrimination. Note also that $\mathcal{D}$ is universally separable in the sense of Pollard [(1984), page 38]. [This follows straightforwardly from Pollard (1984), page 38, problem 4.] Therefore, we conclude that the random variable $\sup_\theta |M_{r,n}(\theta) - M_r(\theta)|$ is measurable. Now, Theorem II.14 from Pollard (1984) implies that

$$
\text{(15)} \qquad \sup_\theta |M_{r,n}(\theta) - M_r(\theta)| = \sup_{D \in \mathcal{D}} |\mathbb{P}_n(D) - \mathbb{P}(D)| \to 0 \qquad \text{a.s.}
$$

Further, write

$$
M_r(\hat{\theta}_{r,n}) - M_r(\theta_0)
$$

$$= M_r(\hat{\theta}_{r,n}) - M_{r,n}(\hat{\theta}_{r,n}) + M_{r,n}(\hat{\theta}_{r,n}) - M_r(\theta_0)$$

$$\le \sup_\theta |M_r(\theta) - M_{r,n}(\theta)| + M_{r,n}(\hat{\theta}_{r,n}) - M_r(\hat{\theta}_{r,n})$$

$$\le 2\sup_\theta |M_r(\theta) - M_{r,n}(\theta)|.$$

Hence, (15) implies

$$|M_r(\hat{\theta}_{r,n}) - M_r(\theta_0)| \to 0, \tag{16}$$

almost surely, as $n \to \infty$. Fix $\varepsilon > 0$. Then by (14) there exists a $\delta > 0$ such that $\max_{\theta \in B_\varepsilon^c(\theta_0)} M_r(\theta) \le M_r(\theta_0) - \delta$. Consequently, we have the set inclusion

$$\{\hat{\theta}_{r,n} \in B_\varepsilon^c(\theta_0) \text{ i.o.}\} \subseteq \{M_r(\hat{\theta}_{r,n}) \le M_r(\theta_0) - \delta \text{ i.o.}\}.$$

But (16) implies that the probability of the event on the right-hand side is zero. Thus, we conclude that $\{\hat{\theta}_{r,n} \in B_\varepsilon(\theta_0) \ ev.\}$ occurs with probability one. Since $\varepsilon > 0$ was arbitrary, we have that $\hat{\theta}_{r,n} \to \theta_0$, almost surely, as $n \to \infty$. This concludes the proof. $\square$

PROOF OF THEOREM 2. The proof is based on verifying conditions of the main theorem of Kim and Pollard (1990) [cf. also Theorem 3.2.10 in van der Vaart and Wellner (1996)].

Let $V(\theta)$ denote the second derivative matrix of the function

$$M_r(\theta) = \mathbb{P}\{|Xa + b - Y| \le r\sqrt{X^2 + 1}\}\mathbb{P}\{|Z^T\theta - Y| \le r\|Z\|\}.$$

Write

$$M_r(\theta) = \mathbb{E}[F(r\|Z\| - Z^T(\theta - \theta_0)) - F(-r\|Z\| - Z^T(\theta - \theta_0))], \tag{17}$$

where $F$ is the distribution function of $\epsilon$, and the expected value above is taken w.r.t. the distribution of $Z := (X, 1)^T$. Now, recall that $f$ is assumed to be continuously differentiable with bounded derivative, and that $\mathbb{E}X^2 < \infty$. Therefore, we can apply the dominated convergence theorem to interchange the order of expectation and differentiation for the expression on the right-hand side of (17). In particular, (17) can be differentiated twice w.r.t. $\theta$ under the integral sign, yielding

$$V(\theta) := \nabla_\theta^2 M_r(\theta)$$
$$= \mathbb{E}\{[f'(r\|Z\| - (\theta - \theta_0)^T Z) - f'(-r\|Z\| - (\theta - \theta_0)^T Z)]ZZ^T\}.$$

Let $V_0 = V(\theta_0)$. Note that the matrix $V_0$ is negative definite when $X$ is nondegenerate. This follows because for strictly unimodal symmetric densities $f$, $f'(x) - f'(-x) < 0$ for all $x > 0$, and under the premise of the theorem, $\mathbb{E}ZZ^T$ is positive definite.

For $\delta > 0$ consider classes of functions $\mathcal{M}_\delta = \{m_\theta - m_{\theta_0} : \|\theta - \theta_0\| < \delta\}$, where $m_\theta = \mathbf{1}_{D_\theta}$, and $D_\theta$ is defined in (4). These classes have polynomial discrimination, that is, they are VC classes [see Pollard (1984), Definition II.13, or van der Vaart and Wellner (1996), page 85] with envelope functions

$$\bar{M}_\delta = \sup_{\|\theta - \theta_0\| < \delta} \left| \mathbf{1}\left\{ -r \leq \frac{Z^T\theta - Y}{\|Z\|} \leq r \right\} - \mathbf{1}\left\{ -r \leq \frac{Z^T\theta_0 - Y}{\|Z\|} \leq r \right\} \right|$$

$$\leq \mathbf{1}\left\{ -r - \delta \leq \frac{Z^T\theta_0 - Y}{\|Z\|} \leq -r + \delta \right\} + \mathbf{1}\left\{ r - \delta \leq \frac{Z^T\theta_0 - Y}{\|Z\|} \leq r + \delta \right\}.$$

Therefore, for small $\delta$,

$$\mathbb{E}\bar{M}_\delta^2 \leq \mathbb{P}\left\{ -r - \delta \leq \frac{\epsilon}{\|Z\|} \leq -r + \delta \right\} + \mathbb{P}\left\{ r - \delta \leq \frac{\epsilon}{\|Z\|} \leq r + \delta \right\}$$

$$\leq c\delta := c\phi^2(\delta)$$

for some positive constant $c$. This verifies condition (vi) in Kim and Pollard [(1990), Theorem 1.1], namely, that $\mathbb{E}\bar{M}_\delta^2 = O(\delta)$. Thus, we anticipate that $n^{-1/3}$ is the rate at which $\hat{\theta}_{r,n}$ converges to $\theta_0$. To arrive at a rigorous conclusion, the key is to compute $\mathbb{E}(m_{\theta_0 + \delta\xi} - m_{\theta_0 + \delta\eta})^2$ for fixed $\delta > 0$ and $\xi, \eta \in \mathbb{R}^2$. This behavior, together with the order of $\phi(\delta)$, will also determine the structure of the increments of the limiting Gaussian process asserted in the theorem. To that end, note that

$$\mathbb{E}[(m_{\theta_0 + \delta\xi} - m_{\theta_0 + \delta\eta})^2 \mathbf{1}\{Z^T\xi \leq Z^T\eta\}]$$

$$= \mathbb{E} \int f(x)\mathbf{1}\{x \in [r\|Z\| + \delta Z^T\xi, r\|Z\| + \delta Z^T\eta]\} \, dx \mathbf{1}\{Z^T\xi \leq Z^T\eta\}$$

$$\quad + \int f(x)\mathbf{1}\{x \in [r\|Z\| + \delta Z^T\xi, r\|Z\| + \delta Z^T\eta]\} \, dx \mathbf{1}\{Z^T\xi \leq Z^T\eta\}$$

$$= \mathbb{E}[F(-r\|Z\| + \delta Z^T\eta) - F(-r\|Z\| + \delta Z^T\xi); Z^T\xi \leq Z^T\eta]$$

$$\quad + \mathbb{E}[F(r\|Z\| + \delta Z^T\eta) - F(r\|Z\| + \delta Z^T\xi); Z^T\xi \leq Z^T\eta]$$

$$=: \mathcal{I}_1 + \mathcal{I}_2.$$

Similar expressions hold when the above expectation is taken on the event $\mathbf{1}\{Z^T\xi > Z^T\eta\}$, with $\xi$ replaced by $\eta$ and vice versa. Our objective is to evaluate an expression for

$$\lim_{\delta \downarrow 0} \frac{\mathbb{E}(m_{\theta_0 + \delta\xi} - m_{\theta_0 + \delta\eta})^2}{\phi^2(\delta)}.$$

But, since $\phi(\delta) = \delta^{1/2}$, this amounts to differentiating $\mathbb{E}(m_{\theta_0 + \delta\xi} - m_{\theta_0 + \delta\eta})^2$ w.r.t. $\delta$ under the integral. (This interchange is justified since $f$, the density of $\epsilon$, is assumed to be bounded, and $Z$ has finite second moment.) Given the above expressions for $\mathcal{I}_1$ and $\mathcal{I}_2$, straightforward algebra yields

$$\lim_{\delta\downarrow 0}\frac{\mathbb{E}(m_{\theta_0+\delta\xi}-m_{\theta_0+\delta\eta})^2}{\phi^2(\delta)}=\mathbb{E}\{[f(-r\|Z\|)+f(r\|Z\|)]|Z^T(\xi-\eta)|\}$$
$$=2\mathbb{E}\{f(r\|Z\|)|Z^T(\xi-\eta)|\}.$$

This completes the proof. □

PROOF OF THEOREM 3. Under the premise of the theorem, there are no parallel lines $L_i$ in the Hough domain. In other words, any pair of random lines intersect, and there is a closed ball of finite radius that contains the set of all intersection points. By construction, for fixed $n$, $\hat{\theta}=\hat{\theta}_{r,n}$ is the center of the ball of radius $r$ that crosses over the maximal number of random lines $L_i$ in the parameter space. Of course, $nM_{r,n}(\hat{\theta}_{r,n})$ is the corresponding number of such lines. Clearly, in order to shift this estimate to infinity one should add at least $nM_{r,n}(\hat{\theta}_{r,n})-1$ lines at infinity. Thus, the smallest contamination fraction under which $\hat{\theta}_{r,n}$ breaks down is $(nM_{r,n}(\hat{\theta}_{r,n})-1)/(n+nM_{r,n}(\hat{\theta}_{r,n})-1)$. Applying the argument as in the proof of Theorem 1, we conclude $M_{r,n}(\hat{\theta}_{r,n})\overset{\text{a.s.}}{\to}M_r(\theta_0)=\mathbb{P}\{\epsilon^2\leq r^2\|Z\|^2\}$, and the result for $\varepsilon_{\text{add}}(\hat{\theta}_{r,n};\mathcal{Y}_n)$ follows. For the replacement breakdown point, it is sufficient to note that under the premise of the theorem at least $\lfloor nM_{r,n}(\hat{\theta}_r,n)/2\rfloor$ lines should be replaced. The proof is complete. □

## APPENDIX B: PROOFS FOR SECTION 5

First we state the uniform central limit theorem for the random field $M_{r,n}(\cdot)$ alluded to before. The statement is formulated in terms of the class of sets generated by the Hough template.

PROPOSITION 2. Let $\mathcal{D}=\{m_\theta=\mathbf{1}_{D_\theta}:\theta\in\mathbb{R}^2\}$, where $D_\theta$ is defined in (4). Let $l^\infty(\mathcal{D})$ denote the set of all uniformly bounded real functions on $\mathcal{D}$. Then the class $\mathcal{D}$ is $\mathbb{P}$-Donsker, that is, $\sqrt{n}(\mathbb{P}_n-\mathbb{P})\Rightarrow G_\mathbb{P}$ in $l^\infty(\mathcal{D})$, where the limit process $\{G_\mathbb{P}m_\theta:m_\theta\in\mathcal{D}\}$ is zero mean Gaussian with covariance function

$$(18)\qquad \mathbb{E}[G_\mathbb{P}m_\xi G_\mathbb{P}m_\eta]=\mathbb{P}(D_\xi\cap D_\eta)-\mathbb{P}(D_\xi)\mathbb{P}(D_\eta).$$

The proposition follows from the uniform central limit theorem for measurable VC-classes [e.g., Corollary 6.3.17 in Dudley (1999)]. Through the mapping $\theta\mapsto D_\theta$, the weak convergence in $l^\infty(\mathcal{D})$ implies that $\sqrt{n}(M_{r,n}(\theta)-M_r(\theta))\Rightarrow G(\cdot)$, where "$\Rightarrow$" denotes weak convergence in $l^\infty(\mathbb{R}^2)$, and the limit is a zero mean Gaussian process with covariance function induced by (18).

PROOF OF THEOREM 4. First we prove the statement given in part (i) of the theorem. The proof proceeds in two steps.

STEP 1. We will require a notion of convergence of sets (all sets are members of the Borel $\sigma$-field over $\mathbb{R}^2$). For any two sets $A_1$, $A_2$, let $A_1 \triangle A_2 := (A_1 \setminus A_2) \cup (A_2 \setminus A_1)$ be the symmetric difference, and define

$$d(A_1, A_2) := \sup_{k \geq 1} \mathcal{L}\{(A_1 \triangle A_2) \cap B_k\},$$

where $\mathcal{L}\{\cdot\}$ stands for Lebesgue measure in $\mathbb{R}^2$, and $B_k = \{\theta \in \mathbb{R}^2 : \|\theta\| \leq k\}$. Note that the above supremum is always finite due to the compactness assumption of the parameter space. First, we prove that

$$(19) \qquad \sup_{\lambda \in \Lambda} d(\Theta_\lambda, \Theta_{\lambda,n}) \to 0 \qquad \text{a.s.,}$$

as $n \to \infty$ [we refer to Molchanov (1998) for closely related results]. For brevity, let us denote $\Delta_{\lambda,n} := \Theta_\lambda \triangle \Theta_{\lambda,n}$. Fix $\delta > 0$. We start with the decomposition

$$d(\Theta_\lambda, \Theta_{\lambda,n}) = \mathcal{L}\{\Delta_{\lambda,n} \cap \{\theta : |M_r(\theta) - \lambda| < \delta\}\}$$
$$+ \mathcal{L}\{\Delta_{\lambda,n} \cap \{\theta : |M_r(\theta) - \lambda| \geq \delta\}\}.$$

The first term on the right-hand side is dominated by $\mathcal{L}\{\{\theta : |M_r(\theta) - \lambda| < \delta\}\}$. The second term on the right-hand side can be upper bounded using the Markov inequality as follows:

$$\mathcal{L}\{\Delta_{\lambda,n} \cap \{\theta : |M_r(\theta) - \lambda| \geq \delta\}\} \leq \delta^{-1} \int_{\Delta_{\lambda,n}} |M_r(\theta) - \lambda| \, d\theta$$
$$\leq \delta^{-1} \mathcal{L}\{\Delta_{\lambda,n}\} \sup_{\theta \in \Delta_{\lambda,n}} |M_r(\theta) - \lambda|.$$

Now, for sufficiently large $n$ (not depending on the choice of $\lambda$) we have (a.s.) the set inclusions

$$\Delta_{\lambda,n} = \{\theta : M_r(\theta) \geq \lambda, M_{r,n}(\theta) < \lambda\} \cup \{\theta : M_r(\theta) < \lambda, M_{r,n}(\theta) \geq \lambda\}$$
$$\subseteq \{\theta : M_r(\theta) \geq \lambda, M_r(\theta) \leq \lambda + \eta_n\} \cup \{\theta : M_r(\theta) < \lambda, M_r(\theta) \geq \lambda - \eta_n\}$$
$$\subseteq \{\theta : |M_r(\theta) - \lambda| \leq \eta_n\},$$

where

$$\eta_n := \sup_{\theta \in \mathbb{R}^2} |M_{r,n}(\theta) - M_r(\theta)|$$

and does not depend on $\lambda$. It follows that

$$\sup_{\theta \in \Delta_{\lambda,n}} |M_r(\theta) - \lambda| \leq \eta_n.$$

In particular, we have for sufficiently large $n$ (independent of $\lambda$) that

$$d(\Theta_\lambda, \Theta_{\lambda,n}) = \mathcal{L}\{\Delta_{\lambda,n}\}$$
$$\leq \mathcal{L}\{\{\theta : |M_r(\theta) - \lambda| \leq \eta_n\}\}$$

and the bound on the right-hand side is uniform in $\lambda$. Thus, taking the supremum over $\lambda \in \Lambda$, letting $n \to \infty$ and appealing to condition (9), we obtain the asserted asymptotic (19).

STEP 2. We now show that for all $\lambda \in \Lambda$,

$$\text{(20)} \qquad \sqrt{n}(E_n(\lambda) - E(\lambda)) = \nu_n(\lambda) + o_p(1), \qquad n \to \infty,$$

where $o_p(1)$ is uniform in $\lambda \in \Lambda$. First, observe that

$$\text{(21)} \qquad \begin{aligned} E_n(\lambda) - E(\lambda) &= \int_{\Theta_{\lambda,n}} (M_{r,n}(\theta) - \lambda)\, d\theta - \int_{\Theta_\lambda} (M_r(\theta) - \lambda)\, d\theta \\ &= \nu_n(\lambda) + R_n, \end{aligned}$$

where

$$\text{(22)} \qquad R_n := \int_{\Theta_{\lambda,n} \setminus \Theta_\lambda} (M_{r,n}(\theta) - \lambda)\, d\theta - \int_{\Theta_\lambda \setminus \Theta_{\lambda,n}} (M_{r,n}(\theta) - \lambda)\, d\theta.$$

Now,

$$\begin{aligned} |\sqrt{n} R_n| &\leq \sqrt{n} \int_{\Theta_\lambda \triangle \Theta_{\lambda,n}} |M_{r,n}(\theta) - \lambda|\, d\theta \\ &\leq d(\Theta_\lambda, \Theta_{\lambda,n}) \sqrt{n} \sup_{\theta \in \Theta_\lambda \triangle \Theta_{\lambda,n}} |M_{r,n}(\theta) - \lambda|. \end{aligned}$$

To prove that $|\sqrt{n} R_n| = o_p(1)$, it suffices to prove this for the right-hand side above. To see this, recall from Step 1 that

$$\sup_{\theta \in \Theta_\lambda \triangle \Theta_{\lambda,n}} |M_r(\theta) - \lambda| \leq \sup_\theta |M_{r,n}(\theta) - M_r(\theta)|,$$

where the upper bound does not depend on $\lambda$. Consequently, we have that

$$|\sqrt{n} R_n| \leq d(\Theta_\lambda, \Theta_{\lambda,n}) \sqrt{n} \sup_\theta |M_{r,n}(\theta) - M_r(\theta)|.$$

But it follows from Proposition 2 that

$$\sqrt{n} \sup_\theta |M_{r,n}(\theta) - M_r(\theta)| \Rightarrow \sup_\theta |G(\theta)|,$$

where $G(\cdot)$ is the zero mean Gaussian process identified in Proposition 2 and the discussion following thereafter, and the above supremum is finite, almost surely. Note that the weak limit does not depend on $\lambda$. By Step 1 we have that $\sup_{\lambda \in \Lambda} d(\Theta_\lambda, \Theta_{\lambda,n}) \to 0$ as $n \to \infty$, a.s. Finally, using Slutzky's lemma, we have that $\sqrt{n}R_n = o_p(1)$ uniformly in $\lambda$. This result, together with (21), gives the assertion (20).

Finally, we put the pieces together using the continuous mapping theorem in the space of continuous functions [see, e.g., Billingsley (1968)], which yields that $\nu_n(\lambda)$ converges to the corresponding integral of the process $G(\cdot)$. To that end, we note that the mapping $\lambda \mapsto \Theta_\lambda$ is continuous w.r.t. the metric $d$, because $M_r(\cdot)$ is continuous and (9) holds. This concludes the proof of the first statement of the theorem.

The proof of statement (ii) goes along the same lines as above. We indicate only the differences. Note that $E_{\mathcal{T}}(\lambda) = E(\lambda)$ because $\Theta_\lambda \in \mathcal{T}$. Also, by definition of $\Theta_{\lambda,n}(\mathcal{T})$,

$$(23) \qquad H_{\lambda,n}\{\Theta_\lambda\} \leq H_{\lambda,n}\{\Theta_{\lambda,n}(\mathcal{T})\} \leq H_{\lambda,n}\{\Theta_{\lambda,n}\}.$$

Therefore, similarly to (21), we write

$$E_{\mathcal{T},n}(\lambda) - E(\lambda) = \nu_n(\lambda) + \tilde{R}_n,$$

where

$$\begin{aligned}
\tilde{R}_n &:= \int_{\Theta_{\lambda,n}(\mathcal{T}) \setminus \Theta_\lambda} (M_{r,n}(\theta) - \lambda)\, d\theta - \int_{\Theta_\lambda \setminus \Theta_{\lambda,n}(\mathcal{T})} (M_{r,n}(\theta) - \lambda)\, d\theta \\
&= H_{\lambda,n}\{\Theta_{\lambda,n}(\mathcal{T})\} - H_{\lambda,n}\{\Theta_\lambda\} \\
&\leq H_{\lambda,n}\{\Theta_\lambda\} - H_{\lambda,n}\{\Theta_{\lambda,n}\} = R_n,
\end{aligned}$$

the last inequality follows from (23) and $R_n$ is defined in (22). Thus, $|\sqrt{n}\tilde{R}_n|$ is bounded using the bounds on $|\sqrt{n}R_n|$ above. Other details of the proof remain unchanged. $\square$

PROOF OF PROPOSITION 1. It follows immediately from the definition that $M_r(\theta) = \tilde{M}_r(\theta - \theta_0)$, where

$$\begin{aligned}
\tilde{M}_r(\theta) &= \mathbb{P}\{|\epsilon + Z^T\theta| \leq r\|Z\|\} \\
&= \mathbb{E}[F(r\|Z\| - Z^T\theta) - F(-r\|Z\| - Z^T\theta)].
\end{aligned}$$

By symmetry of $f$,

$$\begin{aligned}
F(r\|Z\| &- Z^T\theta) - F(-r\|Z\| - Z^T\theta) \\
&= F(r\|Z\| + Z^T\theta) - F(-r\|Z\| + Z^T\theta) \qquad \forall Z,
\end{aligned}$$

and, therefore, $\tilde{M}_r(\theta) = \tilde{M}_r(-\theta)$ $\forall \theta$. Uniqueness of the mode follows from the Anderson lemma.

Let $\theta_1, \theta_2 \in \tilde{\Theta}_\lambda$, that is, $\tilde{M}_r(\theta_1) \geq \lambda$ and $\tilde{M}_r(\theta_2) \geq \lambda$. Let $\theta_* = \alpha\theta_1 + (1-\alpha)\theta_2$ for some $\alpha \in (0,1)$, and denote $I_1 = [-r\|Z\| - Z^T\theta_1, r\|Z\| - Z^T\theta_1]$, $I_2 = [-r\|Z\| - Z^T\theta_2, r\|Z\| - Z^T\theta_2]$, and $I_* = [-r\|Z\| - Z^T\theta_*, r\|Z\| - Z^T\theta_*]$. With this notation,

$$M_r(\theta_*) = \mathbb{E}\int_{I_*} f(x)\,dx.$$

The lengths of $I_1, I_2$ and $I_*$ are equal to $2r\|Z\|$. However, since $\min\{Z^T\theta_1, Z^T\theta_2\} \leq Z^T\theta_* \leq \max\{Z^T\theta_1, Z^T\theta_2\}$, the center of $I_*$ is closer to the origin than one of the centers of $I_1$ and $I_2$. Therefore, by symmetry and unimodality of $f$, for all $Z$,

$$M_r(\theta_*) = \mathbb{E}\int_{I_*} f(x)\,dx \geq \mathbb{E}\min\left\{\int_{I_1} f(x)\,dx, \int_{I_2} f(x)\,dx\right\} \geq \lambda.$$

Thus, $\theta_* \in \tilde{\Theta}_\lambda$, and $\tilde{\Theta}_\lambda$ is convex. $\square$

**Acknowledgments.** The authors wish to thank the referees for their careful reading and helpful and constructive suggestions.

## REFERENCES

ANDERSON, T. W. (1955). The integral of a symmetric unimodal function over a symmetric convex set and some probability inequalities. *Proc. Amer. Math. Soc.* **6** 170–176. MR69229

BILLINGSLEY, P. (1968). *Convergence of Probability Measures.* Wiley, New York. MR233396

CHERNOFF, H. (1964). Estimation of the mode. *Ann. Inst. Statist. Math.* **16** 31–41. MR172382

DANIELS, H. E. (1954). A distribution-free test for regression parameters. *Ann. Math. Statist.* **25** 499–513. MR64374

DAVIES, L. (1990). The asymptotics of $S$-estimators in the linear regression model. *Ann. Statist.* **18** 1651–1675. MR1074428

DONOHO, D. (1988). One-sided inference about functionals of a density. *Ann. Statist.* **16** 1390–1420. MR964930

DONOHO, D. and HUBER, P. J. (1983). The notion of breakdown point. In *A Festschrift for Eric L. Lehmann* (P. J. Bickel, K. A. Doksum and J. L. Hodges, Jr., eds.) 157–184. Wadsworth, Belmont, CA. MR689745

DUDLEY, R. M. (1999). *Uniform Central Limit Theorems.* Cambridge Univ. Press. MR1720712

HARTIGAN, J. A. (1987). Estimation of a convex density contour in two dimensions. *J. Amer. Statist. Assoc.* **82** 267–270. MR883354

HOUGH, P. V. (1959). Machine analysis of bubble chamber pictures. In *International Conference on High Energy Accelerators and Instrumentation* (L. Kowarski, ed.) 554–556. CERN.

GROENEBOOM, P. (1989). Brownian motion with a parabolic drift and Airy functions. *Probab. Theory Related Fields* **81** 79–109. MR981568

GROENEBOOM, P. and WELLNER, J. (2001). Computing Chernoff's distribution. *J. Comput. Graph. Statist.* **10** 388–400. MR1939706

JOHNSTONE, I. M. and VELLEMAN, P. F. (1985). The resistant line and related regression methods (with discussion). *J. Amer. Statist. Assoc.* **80** 1041–1059. MR819612

KIM, J. and POLLARD, D. (1990). Cube root asymptotics. *Ann. Statist.* **18** 191–219. MR1041391

KIRYATI, N. and BRUCKSTEIN, A. M. (1992). What's in a set of points? (Straight line fitting). *IEEE Trans. Pattern Anal. Mach. Intell.* **14** 496–500.

LEAVERS, V. F. (1993). Which Hough transform? *CVGIP: Image Understanding* **58** 250–264.

MOLCHANOV, I. S. (1998). A limit theorem for solutions of inequalities. *Scand. J. Statist.* **25** 235–242. MR1614288

MÜLLER, D. W. and SAWITZKI, G. (1991). Excess mass estimates and tests for multi-modality. *J. Amer. Statist. Assoc.* **86** 738–746. MR1147099

POLLARD, D. (1984). *Convergence of Stochastic Processes.* Springer, New York. MR762984

POLONIK, W. (1995). Measuring mass concentrations and estimating density contour clusters—an excess mass approach. *Ann. Statist.* **23** 855–881. MR1345204

PRINCEN, J., ILLINGWORTH, J. and KITTLER, J. (1992). A formal definition of the Hough transform: Properties and relationships. *J. Math. Imaging Vision* **1** 153–168.

PRINCEN, J., ILLINGWORTH, J. and KITTLER, J. (1994). Hypothesis testing: A framework for analyzing and optimizing Hough transform performance. *IEEE Trans. Pattern Anal. Mach. Intell.* **16** 329–341.

RITTER, G. X. and WILSON, J. N. (1996). *Handbook of Computer Vision Algorithms in Image Algebra.* CRC Press, Boca Raton, FL. MR1398654

ROUSSEEUW, P. J. (1984). Least median of squares regression. *J. Amer. Statist. Assoc.* **79** 871–880. MR770281

ROUSSEEUW, P. J. and HUBERT, M. (1999). Regression depth (with discussion). *J. Amer. Statist. Assoc.* **94** 388–433. MR1702314

ROUSSEEUW, P. J. and LEROY, A. M. (1987). *Robust Regression and Outlier Detection.* Wiley, New York. MR914792

ROUSSEEUW, P. and YOHAI, V. (1984). Robust regression by means of *S*-estimators. In *Robust and Nonlinear Time Series Analysis. Lecture Notes in Statist.* **26** 256–272. Springer, New York. MR786313

RUDIN, W. (1987). *Real and Complex Analysis*, 3rd. ed. McGraw-Hill, New York. MR924157

STEWART, C. V. (1999). Robust parameter estimation in computer vision. *SIAM Rev.* **41** 513–537. MR1719133

VAN DER VAART, A. W. and WELLNER, J. A. (1996). *Weak Convergence and Empirical Processes.* Springer, New York. MR1385671

DEPARTMENT OF STATISTICS
HAIFA UNIVERSITY
HAIFA 31905
ISRAEL
E-MAIL: goldensh@stat.haifa.ac.il

GRADUATE SCHOOL OF BUSINESS
COLUMBIA UNIVERSITY
3022 BROADWAY
NEW YORK, NEW YORK 10027
USA
E-MAIL: assaf@gsb.columbia.edu